\documentclass[journal ]{new-aiaa}
\usepackage[utf8]{inputenc}
\usepackage{textcomp}

\usepackage{graphicx}

\usepackage{amsmath}
\usepackage{amssymb}
\usepackage{amsthm}
\newtheorem{theorem}{Theorem}
\newtheorem{remark}{Remark}
\usepackage[version=4]{mhchem}
\usepackage{siunitx}
\usepackage{longtable,tabularx}
\usepackage{subcaption}
\usepackage{float} 
\renewcommand{\vec}{\boldsymbol}	
\newcommand{\mat}{\boldsymbol}		

\title{Constrained Control Allocation With\\ Continuous-Time Rate Constraints: Three-Dimensional Case}

\author{Süleyman Özkurt \footnote{Research Associate, Institute of Flight Mechanics and Controls, sueleyman.oezkurt@ifr.uni-stuttgart.de}, Adrian Grimm\footnote{Research Associate, Institute of Flight Mechanics and Controls, adrian.grimm@ifr.uni-stuttgart.de} and Walter Fichter\footnote{Professor, Institute of Flight Mechanics and Controls, walter.fichter@ifr.uni-stuttgart.de, AIAA Associate Fellow}}
\affil{University of Stuttgart, 70174 Stuttgart, Germany}
\begin{document}

\maketitle

\section{Introduction}
\lettrine{O}ver-actuated aircraft configurations are particularly prevalent in the eVTOL industry and the military sector, as shown in \cite{Ducard, Durham3}. These configurations have more than three control surfaces to generate three control moments or forces of an aircraft. A control allocation algorithm is required for these aircraft to distribute the commanded control moments of the flight controller to the control surfaces. Several control allocation algorithms exist to distribute the control demand among available actuators. The most common ones are direct allocation (DA) in \cite{Durham1,Durham2,Durham3},  pseudo-inverse (PI) method in \cite{Bordignon,Davidson}, the exact redistributed pseudo-inverse (ERPI) method in \cite{ERP}, the redistributed scaled pseudo-inverse (RSPI) method in \cite{Zhang}, the cascaded generalized inverse (CGI) method in \cite{Baratcart,Page,Acheson}, and solving the control allocation problem, which is an optimization problem, as a quadratic programming (QP) problem in \cite{Oppenheimer, Bodson, Page, Boskovic}. It is expected that these allocation algorithms provide control inputs quickly under the consideration of the actuators' position and rate constraints to ensure a smooth and safe flight of the aircraft. However, current algorithms are not fully capable of covering all of the constraints explicitly.  For a detailed review of all allocation methods the reader is referred to \cite{Durham4,Bodson,Oppenheimer,Johansen}

The DA method, described in \cite[ch.~6.6]{Durham4}, requires the geometric calculation of the attainable moment set (AMS) of the edges, facets, and vertices of the control volume. Afterwards, it is checked for each edge and facet if they intersect with the demanded control moment vector. The computational time of iteratively searching for the intersected edges or facets increases with the number of actuators, and thus, real-time capability cannot always be ensured as discussed in \cite[ch.~6.7]{Durham4} and \cite{Simmons}. Consequently, faster methods to find the correct edges and facets are required.

The PI method is a real-time capable operation; however, it does not consider any constraints. Consequently, calculated control inputs may extend beyond the control volume of the actuators.

The CGI method resolves this issue by eliminating the contributions of saturated actuators from the requested control moments. It performs iterative recalculations of the pseudo-inverse of a reduced control effectiveness matrix. The method is used in various applications such as multicopter systems \cite{Marks,Stepanyan}, ship positioning systems \cite{Shi}, and fighter aircraft \cite{Page,Oppenheimer,Eberhardt} due to its low complexity and fast as well as bounded computation time. Hence, it is favored for safety-critical applications \cite{Raab}.  Compared to the PI method, the CGI can use the entire volume of the AMS as stated in \cite{ERP}. However, the current implementations in \cite{Baratcart,Page} encounter difficulties in preserving the direction of the commanded moments, leading to inappropriate control inputs. A thorough examination of the limitations of the CGI is provided in the reference \cite{ERP}.

 The ERPI method, which is a CGI subvariant, in \cite{ERP} addresses this challenge by transforming the control allocation problem into a direction-preserving, maximum agility allocation problem through the decomposition of the commanded moment into a trim moment and a scalable incremental moment. The solution to this problem (i.e., control inputs) is determined through recursive calculation of the pseudo-inverse. In each recursion step, a saturated actuator is removed from the solution space, and a reduced pseudo-inverse of the control effectiveness matrix is applied to calculate the remaining control inputs. Furthermore, the incremental control moment is scaled recursively to ensure the usage of the entire volume of the AMS and the feasibility of the control inputs. According to \cite{Hafner}, the RSPI method is quite similar to ERPI, and therefore the reader is referred to \cite{Zhang, Hafner} for details.

In the described methods, i.e., DA, PI, CGI, and ERPI, rate constraints are addressed indirectly through the reformulation of these constraints as effective position limits utilizing discretizations as shown in \cite{Oppenheimer, Durham5}. Consequently, allocation methods that consider a continuous and direct formulation of rate constraints are lacking. In such instances, QP offers a substantial advantage due to its capacity to seamlessly incorporate multiple constraints. However, a disadvantage of QP is that a solution can only be calculated if the optimization problem is feasible, i.e., the commanded moments must be attainable under the constraints considered. In order to ensure the feasibility of the QP, constrained slack variables are employed in \cite{Simmons} to transform unattainable moments into attainable ones through scaling. Given that the slack variables must also be solved numerically, the allocation results are predominantly influenced by the empirical calibration of the slack variables' weighting. Consequently, the utilization of the entirety of the AMS is not guaranteed when slack variables are employed; therefore, a diminution in the allocation performance is anticipated. Approaches that exclude slack variables are currently absent.

In summary, established methods such as CGI, ERPI, and DA cannot consider rate constraints in continuous time. QP methods are capable of accomplishing this; however, they are subject to the feasibility condition of the optimization problem, a condition that is presently guaranteed by the introduction of slack variables. To the best of the authors' knowledge, no QP allocation approach has been developed that incorporates continuous-time rate constraints while excluding slack variables. The primary objective of this paper is to propose a design for such a QP allocation approach.\newline This paper offers the following principal contribution: The development of a QP constrained control allocation without slack variables that provide control inputs within the rate and position limits of the actuators. This is achieved by determining the AMS, which is derived from the rate and position constraints, and the point of intersection of the commanded moment vector with the boundary of the AMS. Consequently, the feasibility of the QP is ensured, including both rate and position constraints.

The paper is structured as follows: First, the constrained control allocation problem is reviewed and extended with rate constraints in Sec.~\ref {sec:ConstrainedControlAlloc}. Second, the methodology to determine the AMS and to preserve the direction of control moments in case of unattainable commanded moments is presented in Sec.~\ref {sec:Method}. Third, the suggested methodology is applied in Sec.~\ref{sec:SimResults} to a fighter aircraft control allocation problem of \cite{Durham3}, and it is compared with the ERPI allocation method. Last, conclusions are drawn in Sec.~\ref{sec:Conclusion} based on the results presented in this paper.

\section{Constrained Control Allocation Problem}
\label{sec:ConstrainedControlAlloc}
\begin{figure}[h!]
\centering
\includegraphics[width=.7\textwidth]{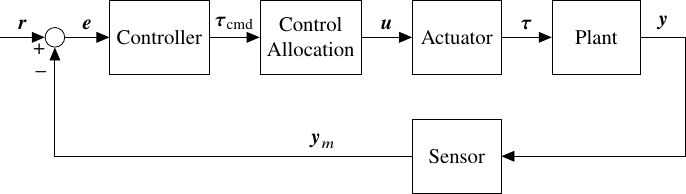}
\caption{Overview of the control architecture.}
\label{fig:AllocIll}
\end{figure}

Within the control allocation problem, the control law, the actuation system, and the physical plant are considered separately as illustrated in Fig.~\ref{fig:AllocIll}.
The controller calculates the demanded control moment $\vec{\tau}_{\mathrm{cmd}}$ based on the measured states $\vec{y}_{m}$ and the commanded reference values $\vec{r}$. The control allocation distributes the demanded control moments to the actuators and generates actuator commands $\vec{u}$. These actuator commands are linearly mapped\footnote{Under the hypothesis that actuator dynamics are significantly faster than system dynamics, a static and linear actuation is considered. Further, $\vec{u}_{\mathrm{act}}=\vec{u}$ is assumed in Fig. \ref{fig:AllocIll}, indicating that the actuators directly implement the allocated control inputs.} onto the pseudo control moment $\vec{\tau}$ with
\begin{equation}
    \vec{\tau}=\mat{B}\vec{u},
    \label{eq:LinMap}
\end{equation}
where $\mat{B} \in \mathbb{R}^{o\times m}$ is the control effectiveness matrix. The nonlinear plant dynamics are given as 
\begin{equation}
\begin{split}
\vec{\dot{x}}&=\vec{f}\left(\vec{x},\vec{\tau}\right),\\
\vec{y}&=\vec{h}\left(\vec{x}\right),
\end{split}
    \label{eq:Dynamics}
\end{equation} 
where $\vec{f}:\mathbb{R}^{n}\times\mathbb{R}^{o}\rightarrow\mathbb{R}^{n}$ is a nonlinear function that maps the state vector $\vec{x}$ and pseudo-control vector $\vec{\tau}$ (e.g., forces and torques for an aircraft) to the state vector derivative. $\vec{h}:\mathbb{R}^{n}\rightarrow\mathbb{R}^{p}$ is the output function and $\vec{y} \in \mathbb{R}^{p}$ is the output vector, respectively.
In practical actuator systems, the actuator command $\vec{u}$ is subject to position constraints. The feasible set $\mat{U}$ for the position limits is defined as
\begin{equation}
    \mat{U}:=\left\{\vec{u}\in\mathbb{R}^{m}\big|\mat{M}\vec{u}\preceq\vec{b}_{u}\right\}\subset \mathbb{R}^{m},
    \label{eq:FeasU}
\end{equation}
where the position constraints are represented as a matrix inequality. The symbol $\preceq$ is employed to denote row-wise inequality. $\mat{U}$ corresponds to the control volume of the actuator. 
\begin{remark}
The feasible set $\mat{U}$ is represented by multiple intersections of half-spaces and hyperplanes (which are convex), providing a compact and bounded set in $\mathbb{R}^{m}$. Hence, it is considered a convex polytope, see  \cite[ch.~2-3]{Grünbaum}.
\end{remark}
The attainable moment set $\mat{D}_{\tau_{u}}$ is defined as 
\begin{equation}
    \mat{D}_{\tau_{u}}:=\left\{\vec{\tau}\in\mathbb{R}^{o}\big|\vec{\tau}=\mat{B}\vec{u}~\land~\vec{u} \in \mat{U} \right\} \subset \mathbb{R}^{o}
\end{equation}
\begin{theorem}
 The attainable moment set $\mat{D}_{\tau_{u}}$ is a convex polytope.    
\end{theorem}
\begin{proof}
Since $\vec{u}$ is limited within the convex polytope $\mat{U}$, the linear transformation $\vec{\tau}=\mat{B}\vec{u}$ will preserve the convexity and will map to a new polytope, see \cite[p.~175]{Köthe}.
\end{proof}
For over-actuated systems, i.e., $m>o$, Eq. \eqref{eq:LinMap} cannot be solved directly. Hence, an actuator command $\vec{u}$ is found by solving the equation below
\begin{equation}
\begin{split}
\mathrm{min}&~\vec{u}^{\intercal}\mat{R}\vec{u}\\
\mathrm{s.t}&~ \vec{\tau}=\mat{B}\vec{u},\\ 
&~ \mat{M}\vec{u}\preceq\vec{b}_{u},
\end{split}
\label{eq:ConstrainedControlAlloc}
\end{equation}
where $\mat{R}$ is a positive definite weighting matrix. Throughout this document, $\mat{R}$ corresponds to the identity matrix. 

The optimization problem in Eq.~\eqref{eq:ConstrainedControlAlloc} is known as the constrained control allocation problem.
\subsection{Constrained Control Allocation With Rate Constraints}
For real actuation systems, the dynamics of the actuator must be considered, which means that a real actuator may fail to implement the commanded control input $\vec{u}$ due to inappropriate control input rates $\vec{\dot{u}}$.
 Respectively, the feasible set $\mat{\dot{U}}$ for the rate limits is defined as
\begin{equation}
\mat{\dot{U}}:=\left\{\vec{\dot{u}}\in\mathbb{R}^{m}\big|\mat{Z}\vec{\dot{u}}\preceq\vec{b}_{\dot{u}}\right\}\subset \mathbb{R}^{m}
\end{equation}
\begin{remark}
The feasible set $\mat{\dot{U}}$ is represented by multiple intersections of half-spaces and hyperplanes (which are convex), providing a compact and bounded set in $\mathbb{R}^{m}$. Hence, it is considered a convex polytope, see  \cite[ch.~2-3]{Grünbaum}.
\end{remark}
To include rate constraints into the constrained control allocation of Eq.~\eqref{eq:ConstrainedControlAlloc}, a continuous first-order relationship between $\vec{u}$ and $\dot{\vec{u}}$ is assumed.  Based on this assumption,  Eq.~\eqref{eq:ConstrainedControlAlloc} is augmented as follows
\begin{equation}
    \begin{split}
        \mathrm{min}~& \begin{bmatrix}
            \vec{u}\\
            \vec{\dot{u}}
        \end{bmatrix}^{\intercal} \mat{R} \begin{bmatrix}
        \vec{u}\\
        \vec{\dot{u}}\\
        \end{bmatrix}\\
        \mathrm{s.t.}~& \vec{\tau}=\mat{B}\vec{u},\\
          \vec{\dot{u}}&=\mat{A}\vec{u},\\
        ~& \mat{M}\vec{u}\preceq \vec{b}_{u},\\
        ~& \mat{Z}\vec{\dot{u}}\preceq \vec{b}_{\dot{u}},\\
    \end{split}
    \label{eq:ConstrainedAlloc}
\end{equation}
where $\vec{\dot{u}}=\mat{A}\vec{u}$ induces a homogeneous linear relationship of the first order with $\mat{A} \in \mathbb{R}^{m\times m}$. This dynamic is imposed on the optimization problem. 

Under the assumption that $\mathrm{det}(A)\neq 0$, the attainable moment set $\mat{D}_{\tau_{\dot{u}}}$ is restricted due to the feasible set $\mat{\dot{U}}$ as defined below
\begin{equation}
    \mat{D}_{\tau_{\dot{u}}}:=\left\{\vec{\tau}\in\mathbb{R}^{o}\big|\vec{\tau}=\mat{B}\mat{A}^{-1}\vec{\dot{u}}~\land~\vec{\dot{u}} \in \mat{\dot{U}} \right\} \subset \mathbb{R}^{o}
\end{equation}
\begin{theorem}
 The attainable moment set $\mat{D}_{\tau_{\dot{u}}}$ is a convex polytope.    
\end{theorem}
\begin{proof}
Since $\vec{\dot{u}}$ is limited within the convex polytope $\mat{\dot{U}}$, the linear transformation $\vec{\tau}=\mat{B}\mat{A}^{-1}\vec{\dot{u}}$ will preserve the convexity and will map to a new polytope, see \cite[p.~175]{Köthe}
\end{proof}
Similar to the discrete approaches in \cite{Oppenheimer}, the additional dynamical constraint leads to new position limits induced by actuator rate limits. However, this formulation does not require a measurement of the current actuator position, which is a major benefit. To guarantee a solution to the optimization problem in Eq. \eqref{eq:ConstrainedAlloc}, the feasibility of the problem must be ensured. This means that the commanded moments $\vec{\tau}_{\mathrm{cmd}}$ must be an element of the attainable moment set $\vec{\tau} \in \mat{D}_{\tau}$ with $\mat{D}_{\tau}:=\left\{\mat{D}_{\tau_{u}} \cap \mat{D}_{\tau_{\dot{u}}}\right\}$\footnote{$\mat{D}_{\tau}$ is a convex polytope, see \cite[ch.~2.3]{Boyd}.}. In case of unattainable commanded moments, this is achieved by the direction-preserving scaling method described below. 
\section{Methodology}
\label{sec:Method}
In this section, the calculation of the AMS, based on rate and position constraints, is presented. Additionally, in Sec.~\ref{sec:MethodA}, a method to preserve the direction of unattainable commanded moments for the three-dimensional case is described.

The calculation of the attainable moment set $\mat{D}_{\tau}$ is based on the calculation of the convex hulls $\mat{D}_{\tau_{u}}\subseteq\tilde{\mat{D}}_{\tau_{u}}$ and $\mat{D}_{\tau_{\dot{u}}}\subseteq\tilde{\mat{D}}_{\tau_{\dot{u}}}$. 

To calculate the convex hull $\tilde{\mat{D}}_{\tau_{u}}$, the vertices of the control volume $\mat{U}$ are used. Vertices of the control volume contain maximum and minimum values of the control inputs $\vec{u}$. The set of vertices is defined as 
\begin{equation}
\mat{V}(\mat{U}):=\left\{\vec{v}_{1},...,\vec{v}_{2^{m}}\right\},
\end{equation}
where $\vec{v}_{i}=[v_{i1},...,v_{im}]^{\intercal}; i\in [1,...,2^{m}]$ and $v_{ij}\in[\mathrm{max}(u_j);\mathrm{min}(u_j)];j \in [1,...,m]$.

The vertices represent control inputs at the boundary of the control volume $\partial(\mat{U})$. As shown in \cite[ch.~5.2]{Durham4}, not all vertices map to the boundary of the attainable moment set $\partial(\mat{D}_{\tau_{u}})$. Some of the vertices will map into the interior of $\mat{D}_{\tau_{u}}$, while a few will map to points on $\partial(\mat{D}_{\tau_{u}})$. This property is used below to calculate $\tilde{\mat{D}}_{\tau_{u}}$, since the convex hull contains the boundary $\partial(\mat{D}_{\tau_{u}})$,
\begin{equation}
\begin{split}
\mat{V}_{D}:=&\left\{\vec{\tau}\in\mathbb{R}^{o}\big|\vec{\tau}=\mat{B}\vec{v} \land \vec{v}\in \mat{V}\right\},\\
\tilde{\mat{D}}_{\tau_{u}}=&~\mathrm{Conv}(\mat{V}_{D}),
\end{split}
\end{equation}
where $\mat{V}_{D}$ is the set of attainable moments on the $\partial(\mat{D}_{\tau_{u}})$ and the interior of $\mat{D}_{\tau_{u}}$. $\mathrm{Conv}(\cdot)$ is the convex hull operator.
The same procedure is applied to calculate $\tilde{\mat{D}}_{\tau_{\dot{u}}}$ as shown below
\begin{equation}
\begin{split}
    \mat{Z}(\mat{\dot{U}}):=&\{\vec{z}_{1},...,\vec{z}_{2^{m}}\},\\
    \mat{Z}_D:=&\left\{\vec{\tau}\in\mathbb{R}^{o}\big|\vec{\tau}=\mat{B}\mat{A}^{-1}\vec{z}\land \vec{z}\in\mat{Z}\right\},\\
    \tilde{\mat{D}}_{\tau_{\dot{u}}}=&~\mathrm{Conv}(\mat{Z}_D),
\end{split}
\end{equation}
where $\mat{Z}(\mat{\dot{U}})$ is the set of vertices of the feasible set $\mat{\dot{U}}$ with $\vec{z}_{i}=[z_{i1},...,z_{im}]^{\intercal}; i\in [1,...,2^{m}]$ and $z_{ij}\in[\mathrm{max}(\dot{u}_j);\mathrm{min}(\dot{u}_j)];j \in [1,...,m]$.
Afterwards, the AMS is obtained through the intersection of both convex hulls $\mat{D}_{\tau}:={\tilde{\mat{D}}_{\tau_{u}}\cap\tilde{\mat{D}}_{\tau_{\dot{u}}}}$. For the calculation of the convex hull, several algorithms are available as discussed in \cite{Qhull,Nielsen}.

\subsection{Direction Preservation of Commanded Control Moments for the Three-Dimensional Case}
\label{sec:MethodA}
\begin{figure}[h!]
\centering
\includegraphics[width=.4\textwidth]{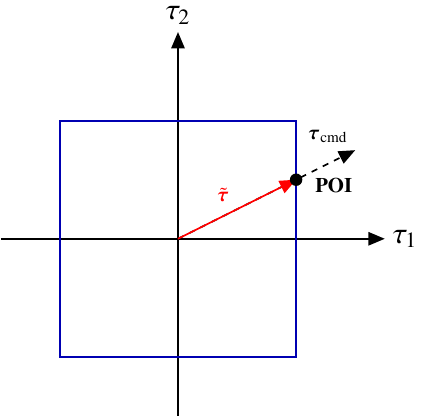}
\caption{Two-dimensional example of the direction-preserving method of the commanded vector.}
\label{fig:AMSDraw}
\end{figure}

The commanded moments $\vec{\tau}_{\mathrm{cmd}}$ must be feasible to obtain a solution for the control inputs $\vec{u}$ from the optimization in Eq.~\eqref{eq:ConstrainedAlloc}.  The feasibility is ensured through the scaling of unattainable moments to the boundary of the AMS. In the three-dimensional case, this is performed by calculating the point of intersection (POI) of the commanded control moment vector with the facets of the AMS as illustrated in Fig.~\ref{fig:AMSDraw}.

 As shown in Fig.~\ref{fig:AMS}, $\vec{\tau}_{\mathrm{cmd}}$ intersects the AMS's facet (blue square) at the POI. If the POI is known, the commanded moment can be clipped to $\tilde{\vec{\tau}}$. To calculate the POI, the reconstructed AMS is rewritten in half-space form (i.e., as a matrix inequality) representing a convex polytope of the form
 \begin{equation}
     \mat{N}\vec{\tau}\preceq\vec{1},
\label{eq:AMSDesc}
 \end{equation}
 where $\vec{1}$ is vector containing only ones and $\mat{N}\in \mathbb{R}^{k\times o}$. Note that in scenarios involving more than three dimensions, the determination of $\mat{N}$ can pose significant computational challenges as discussed in \cite{Facet}. 
If $\vec{\tau}_{\mathrm{cmd}}$ is outside the polytope, then $\mat{N}\vec{\tau}_{\mathrm{cmd}}\npreceq\vec{1}$. Hence, a scaling vector $\vec{r}_{s}$, containing a scaling factor for each inequality, is calculated below.
\begin{equation}
\vec{r}_{s}=\vec{1}\oslash(\mat{N}\vec{\tau}_{\mathrm{cmd}}),
\end{equation}
where $\oslash$ is the operator for element-wise division. Note that $\vec{r}$ represents the distances of $k$ hyperplanes from the origin.
The scalar vector is used to calculate $\mat{\Omega} \in \mathbb{R}^{o\times k}$ below
\begin{equation}
\mat{\Omega}=\left\{\vec{\tilde{\zeta}}_{1},...,\vec{\tilde{\zeta}}_{k}\right\}=\vec{\tau}_{\mathrm{cmd}}\odot\vec{r}_{s}^{\intercal}
\end{equation}
where $\odot$ represents the element-wise multiplication operator. $\mat{\Omega}$ contains moment vectors that are multiplied by $\mat{N}$ to check which scaling factors lead to control moments that comply with all inequalities (i.e., obtaining control moments within or at the boundary of the AMS). This can result in several points of intersection with the boundary of the AMS \footnote{As is obvious from Fig.~\ref{fig:AMS}, there can be multiple POIs with the AMS, if $\vec{\tau}_{\mathrm{cmd}}$ is considered as the directional vector of a straight line through the origin.}.  Hence, the obtained set of POIs is defined as  
\begin{equation}
    \mathrm{POIs}:=\left\{\vec{\tilde{\zeta}_{i}}\in\mathbb{R}^{o}\big|\mat{N}\vec{\tilde{\zeta}}_i\preceq\vec{1}; i\in [1,..,k]\right\}=\left\{\vec{\tilde{\tau}_{1}},...,\vec{\tilde{\tau}}_{n}\right\}.
\end{equation}
The correct POI for the preservation of the direction is found by selecting the moment vector as described below
\begin{equation}
\mathrm{POI}=\vec{\tau}=\underset{\vec{\tilde{\tau} \in \mathrm{POIs}}}{\arg\min}||\vec{\tau}_{\mathrm{cmd}}-\vec{\tilde{\tau}}||_{2}.
\label{eq:POI}
\end{equation}
If $\vec{\tau}_{\mathrm{cmd}}$ is within the polytope, i.e., $\mat{N}\vec{\tau}_{\mathrm{cmd}}\preceq\vec{1}$ then 
\begin{equation}
    \vec{\tau}=\vec{\tau}_{\mathrm{cmd}}.
    \label{eq:tau}
\end{equation}
The moment vector determined from Eq.~\eqref{eq:POI} or Eq.~\eqref{eq:tau} is then passed to the optimization problem in Eq.~\eqref{eq:ConstrainedAlloc}, which can be solved using a QP solver such as \textit{quadprog} \cite{quadprog} or \textit{qpoases} \cite{qpoases1,qpoases2}. 
\section{Simulation Results}
\label{sec:SimResults}
In this section, the proposed approach of Sec.~\ref{sec:Method} is applied to an F18 fighter jet from \cite{Durham1} as an empirical case study to analyze its performance. Further, in Sec.~\ref{sec:SimRes1} the results of the proposed approach are compared with the ERPI approach \cite{ERP}. In Sec.~\ref{sec:SimRes2}, the allocation results provided by the proposed approach are analyzed. These results are considered in the context of real actuator dynamics, including position and rate constraints.

\begin{figure}[h!]
    \centering
    \begin{subfigure}[b]{0.4\textwidth}
        \includegraphics[width=\linewidth]{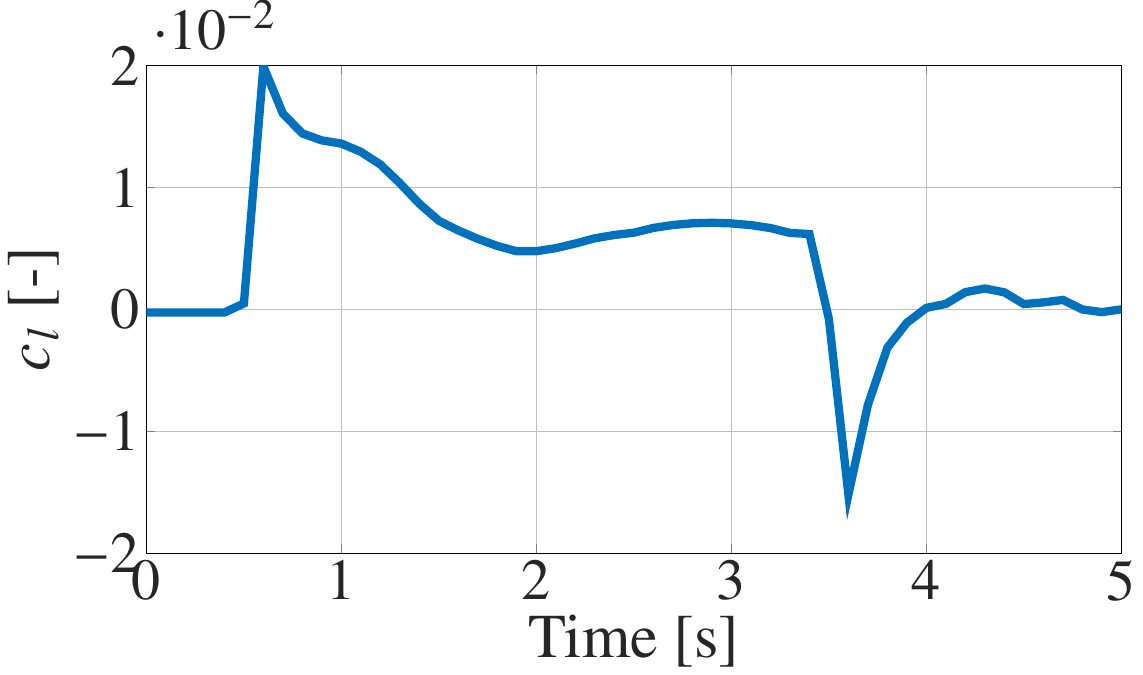}
        \caption{Commanded roll moment.}
        \label{fig:cl}
    \end{subfigure}
        \hfill 
    \begin{subfigure}[b]{0.4\textwidth}
        \includegraphics[width=\linewidth]{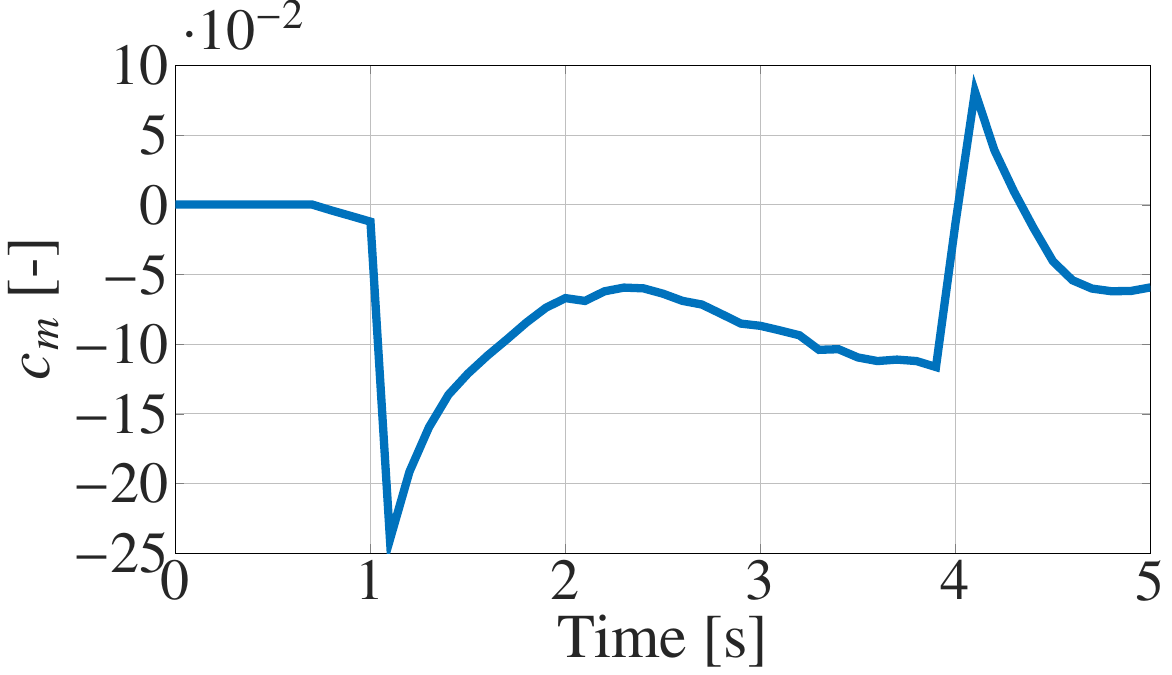}
        \caption{Commanded pitch moment.}
        \label{fig:cm}
    \end{subfigure}
        \hfill
    \begin{subfigure}[b]{0.4\textwidth}
        \includegraphics[width=\linewidth]{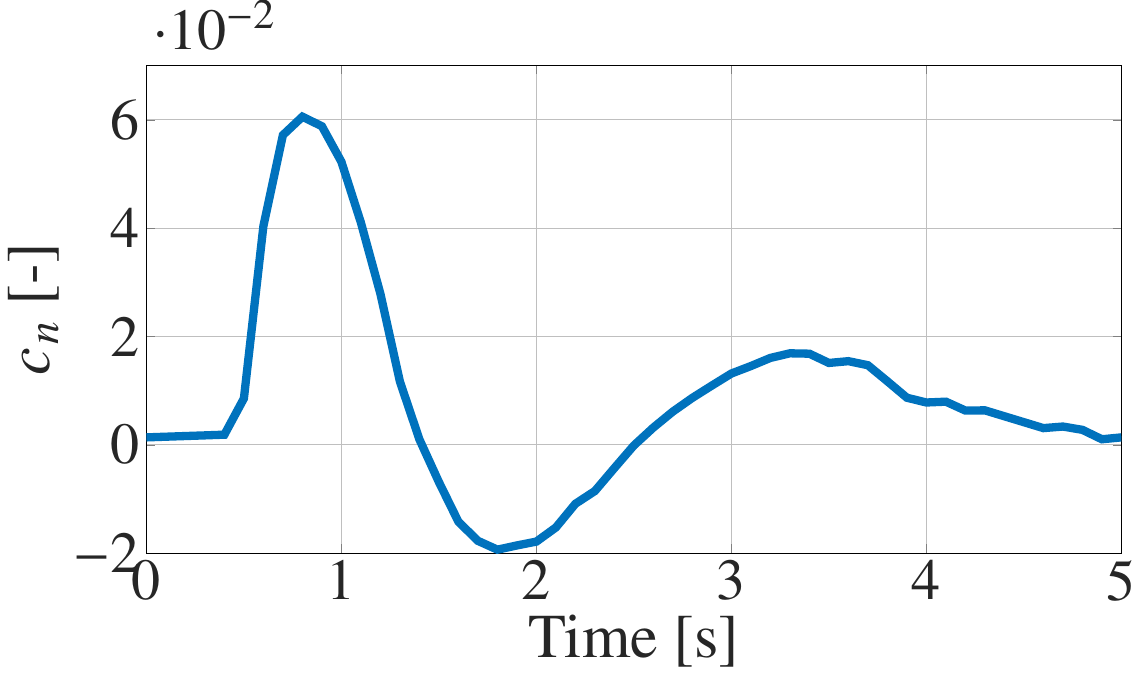}
        \caption{Commanded yaw moment.}
        \label{fig:cn}
    \end{subfigure}
    \caption{Commanded maneuver illustrated as coefficients for the moments.}
    \label{fig:Maneuver}
\end{figure}

The control effectiveness matrix of the F18 aircraft is given as
\begin{equation}
\mat{B}=\frac{10^{-5}}{\mathrm{deg}}\cdot
\begin{array}{c}
\begin{array}{rrrrrrr}
[\mathrm{tail_{l}} & \mathrm{tail_{r}} & \mathrm{flap_{l}} &\mathrm{flap_{r}} & \mathrm{ail_{l}} &\mathrm{ail_{r}} & \mathrm{rudder}]
\end{array}\\\overbrace{
\left[
\begin{array}{rrrrrrr}
23.8 & -23.8 & 123.0 & -123.0 & 41.8 & -41.8 & 3.6 \\
-698.0 & -698.0 & 99.4 & 99.4 & -55.2 & -55.2 & 0.0 \\
-30.9 & 30.9 & 0.0 & 0.0 & -17.4 & 17.4 & -56.2
\end{array}
\right]}
\end{array}
\label{eq:B}
\end{equation}
According to \cite{Harv}, the upper $\vec{\bar{u}}$ and lower $\vec{\underline{u}}$ position limits of the actuator are given as
\begin{equation}
    \begin{split}
        \vec{\bar{u}}&=\begin{bmatrix}
            10.5 & 10.5&45.0&45.0&42&42&30
        \end{bmatrix}^{\intercal}~\mathrm{deg},\\
        \vec{\underline{u}}&=\begin{bmatrix}
            -24&-24&-8&-8&-25&-25&-30
        \end{bmatrix}^{\intercal}~\mathrm{deg}.\\
    \end{split}
\end{equation}
The upper $\vec{\dot{\bar{u}}}$ and lower $\vec{\dot{\underline{u}}}$ rate limits of the actuators are given as 
\begin{equation}
    \begin{split}
        \vec{\dot{\bar{u}}}&=\begin{bmatrix}
            40 & 40&18&18&100&100&82
        \end{bmatrix}^{\intercal}~\frac{\mathrm{deg}}{\mathrm{s}},\\
        \vec{\dot{\underline{u}}}&=\begin{bmatrix}
            -40&-40&-18&-18&-100&-100&-82
        \end{bmatrix}^{\intercal}~\frac{\mathrm{deg}}{\mathrm{s}}.\\
    \end{split}
\end{equation}
The actuator dynamics of the control surfaces are represented by second-order systems
\begin{equation}
  \mathrm{Act}(s)=\frac{\omega_{0}^2}{s^2+2\zeta\omega_{0}s+\omega_{0}^2},
\end{equation} 
with the parameter values below 
\begin{table}[h!]
\centering
\caption{Actuator parameters of the control surfaces.}
\begin{tabular}{|p{3cm}|p{3cm}|p{3cm}|p{3cm}|p{3cm}|}
\hline 
Parameter\textbackslash Surface &Tail stabilizer & Rudder & Aileron & Flaps \\
\hline 
$\omega_{0}$&30.74&72.1&75&35\\
\hline
$\zeta$&0.509&0.69&0.59&0.71\\ \hline
\end{tabular}
\label{tab:noten}
\end{table}
\newline A five-second sequence of the commanded moment vector from \cite{Durham3}, shown in Fig.~\ref{fig:Maneuver}, is used as input for the control allocation problem. The moment vector $\vec{\tau}=[c_{l},c_{m},c_{n}]^{\intercal}$ consists of the roll moment $c_{l}$, the pitch moment $c_{m}$, and the yaw moment $c_{n}$ coefficients. The sequence begins with a $5~g$ level turn to the left. An aggressive roll maneuver to the right is performed within $0.5~\mathrm{s}<t<3.5~\mathrm{s}$. During the same time, around $1~\mathrm{s}<t<4~\mathrm{s}$, the load factor decreases to $3~g$ before it increases to $5~g$ again. The maneuver ends with a level turn to the right after $t\geq 3.5~\mathrm{s}$.
\subsection{Comparison of Methodologies}
\label{sec:SimRes1}
In this section, the outcomes and performance of the proposed approach of Sec.~\ref{sec:Method} are compared with the ERPI algorithm from \cite{ERP} since it is field-proven in several flight tests \cite{Stephan,Ole}, such as the RSPI algorithm in \cite{Surman}.
For the comparison, only the position constraints of the actuator are considered, corresponding to the control allocation problem of Eq. \eqref{eq:ConstrainedControlAlloc}. In addition, actuator dynamics are neglected, assuming that they directly implement the calculated control inputs from the allocation problem, i.e., $\vec{u}=\vec{u}_{\mathrm{act}}$. 

\begin{figure}[H]
\centering
\includegraphics[width=0.6\textwidth]{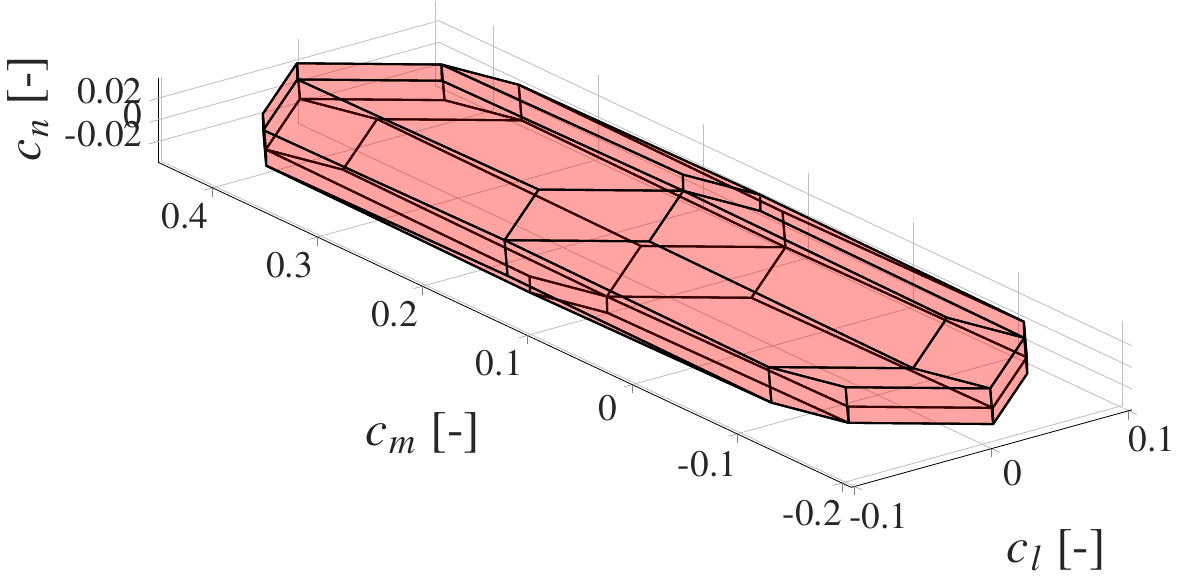}
\caption{Determined convex hull $\tilde{D}_{\tau_{u}}$ from the vertices of $\mat{U}$.}
\label{fig:AMS}
\end{figure}

The suggested methodology of Sec.~\ref{sec:Method} provides the convex hull $\tilde{D}_{\tau_{u}}$ in Fig.~\ref{fig:AMS}. The calculation of the convex hull is performed using the \textit{qhull} algorithm from \cite{Qhull} in MATLAB. The determined $\tilde{D}_{\tau_{u}}$ corresponds to the true AMS of \cite{Durham4}, indicating that the method is capable of adequately reproducing it.  

\begin{figure}[H]
    \centering
    \begin{subfigure}[b]{0.45\textwidth}
        \includegraphics[width=\textwidth]{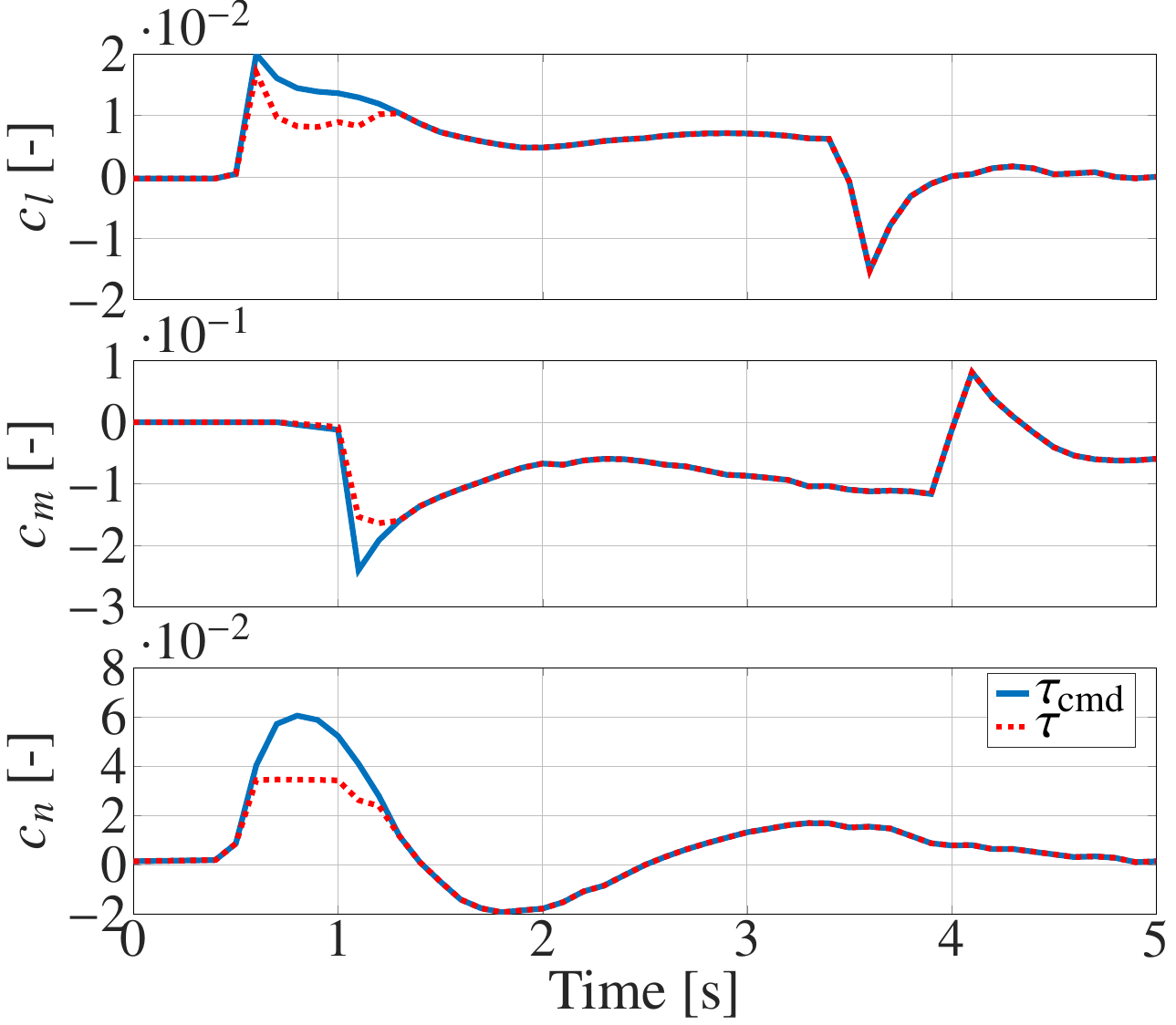}
        \caption{Attainable moments provided by the ERPI algorithm.}
        \label{fig:ERPAM}
    \end{subfigure}
    \hfill
    \begin{subfigure}[b]{0.45\textwidth}
        \includegraphics[width=\textwidth]{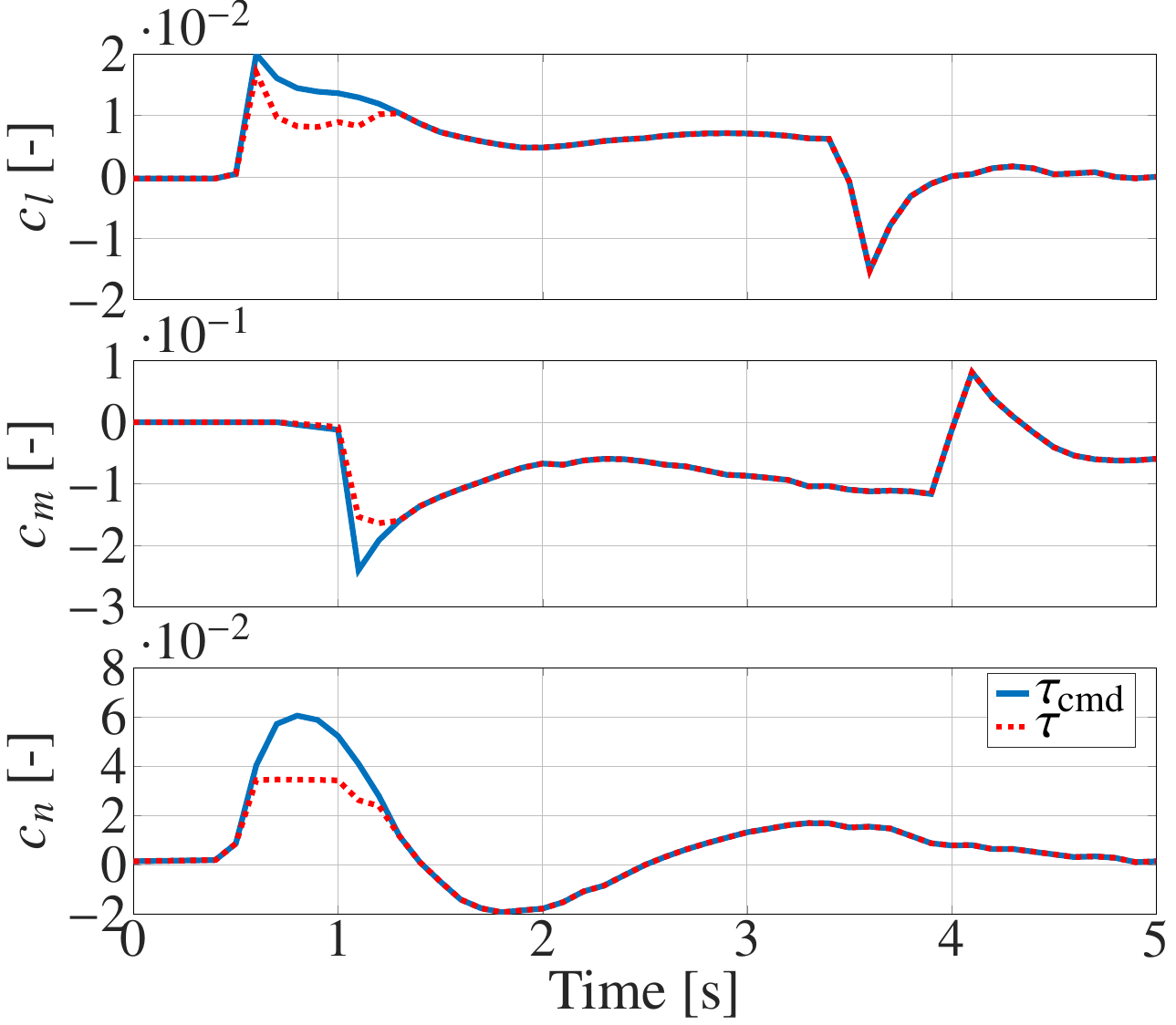}
        \caption{Attainable moments provided by our approach.}
        \label{fig:QPAM}
    \end{subfigure}
    \caption{Comparison of the achieved attainable moments.}
    \label{fig:AMComp}
\end{figure}

The \textit{qpoases} algorithm from \cite{qpoases1,qpoases2} is used for the numerical solving of the QP formulated in Eq.~\eqref{eq:ConstrainedControlAlloc}. 
The determined control moments from both allocation approaches are illustrated in Fig.~\ref{fig:AMComp}.  Since the ERPI algorithm also preserves the direction of the unattainable control moments, the results in Fig.~\ref{fig:AMComp} indicate that our methodology can also preserve the direction of unattainable control moments.

\begin{figure}[H]
    \centering
    \begin{subfigure}[b]{0.45\textwidth}
        \includegraphics[width=\textwidth]{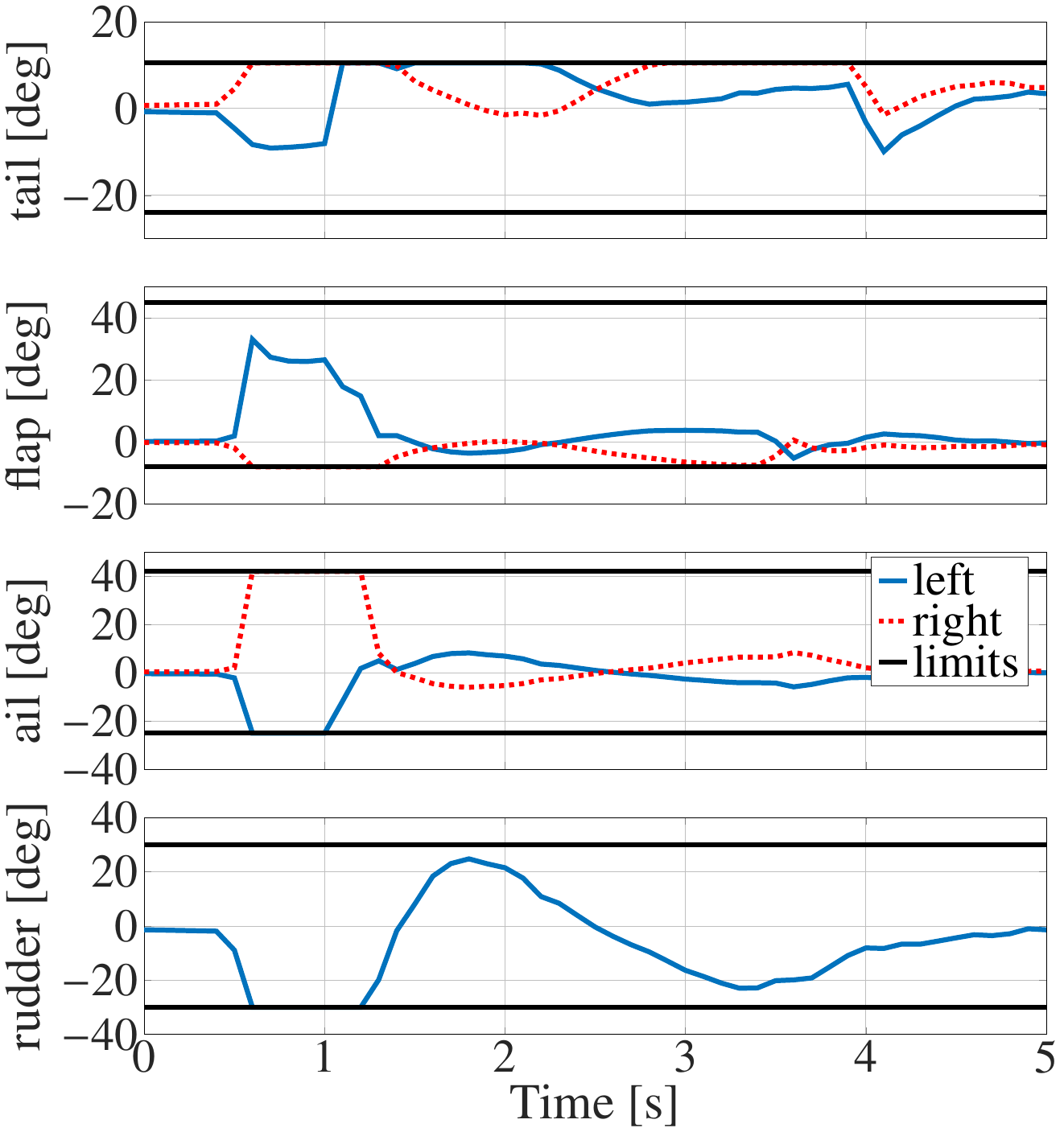}
        \caption{Control inputs provided by the ERPI allocation.}
        \label{fig:ERPU}
    \end{subfigure}
    \hfill
    \begin{subfigure}[b]{0.45\textwidth}
        \includegraphics[width=\textwidth]{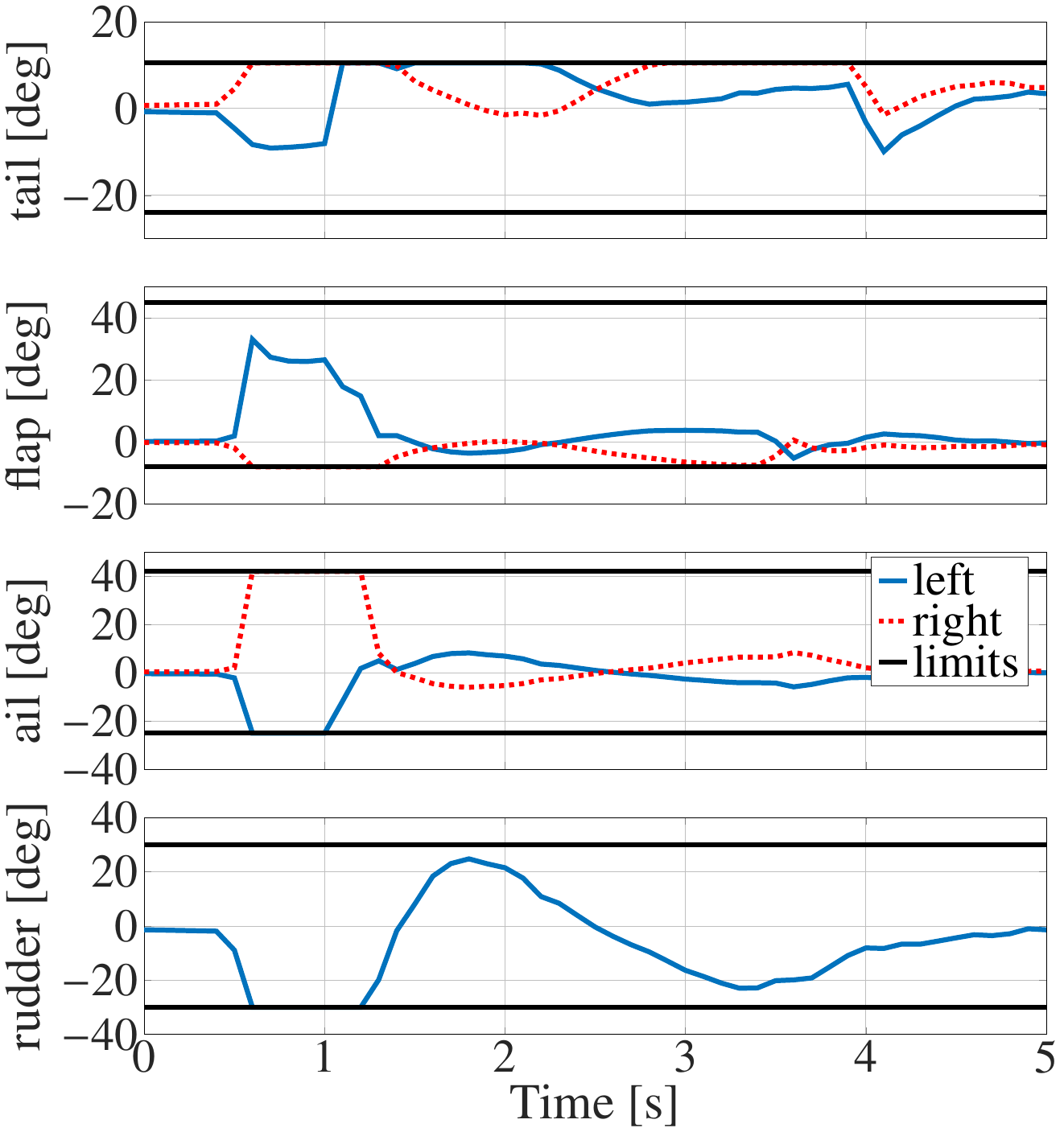}
        \caption{Control inputs provided by our approach.}
        \label{fig:QPU}
    \end{subfigure}
    \caption{Comparison of the calculated control inputs.}
    \label{fig:UComp}
\end{figure}

 Furthermore, since the calculated achievable moments are equal for both methods, the determined control inputs in Fig.~\ref{fig:UComp} are also equal for both methods. Both methods provide control inputs within the position limits of the actuators. Based on the outcomes in Figs.~\ref{fig:AMComp} - \ref{fig:UComp}, the proposed methodology yields results equivalent to the contemporary ERPI method.

\begin{figure}[H]
\centering
\begin{subfigure}[b]{0.45\textwidth}
\includegraphics[width=\textwidth]{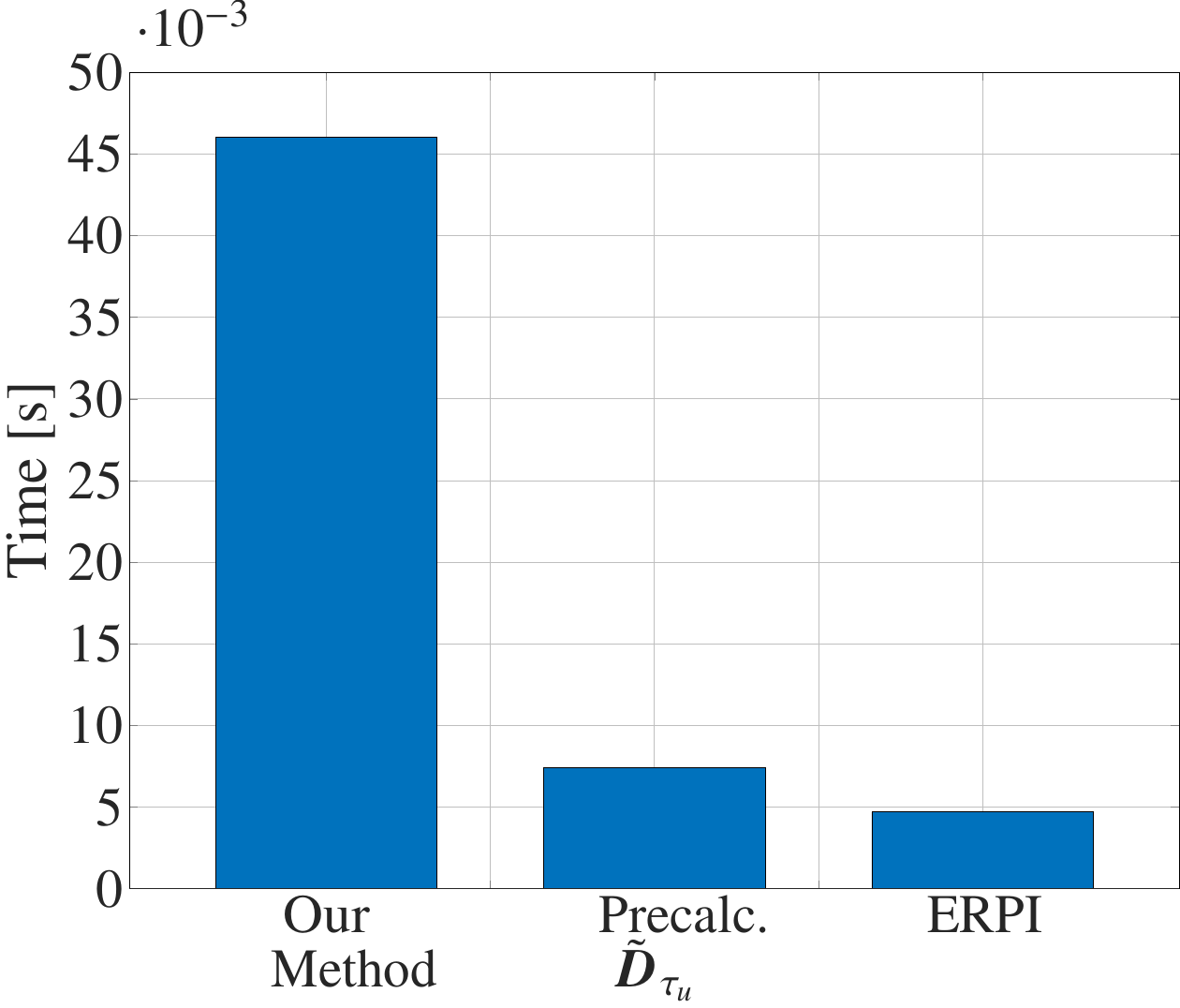}
\caption{Average elapsed time over 100 runs of the allocation algorithms.}
\label{fig:Time}
\end{subfigure}
\hfill
\begin{subfigure}[b]{0.45\textwidth}
\includegraphics[width=\textwidth]{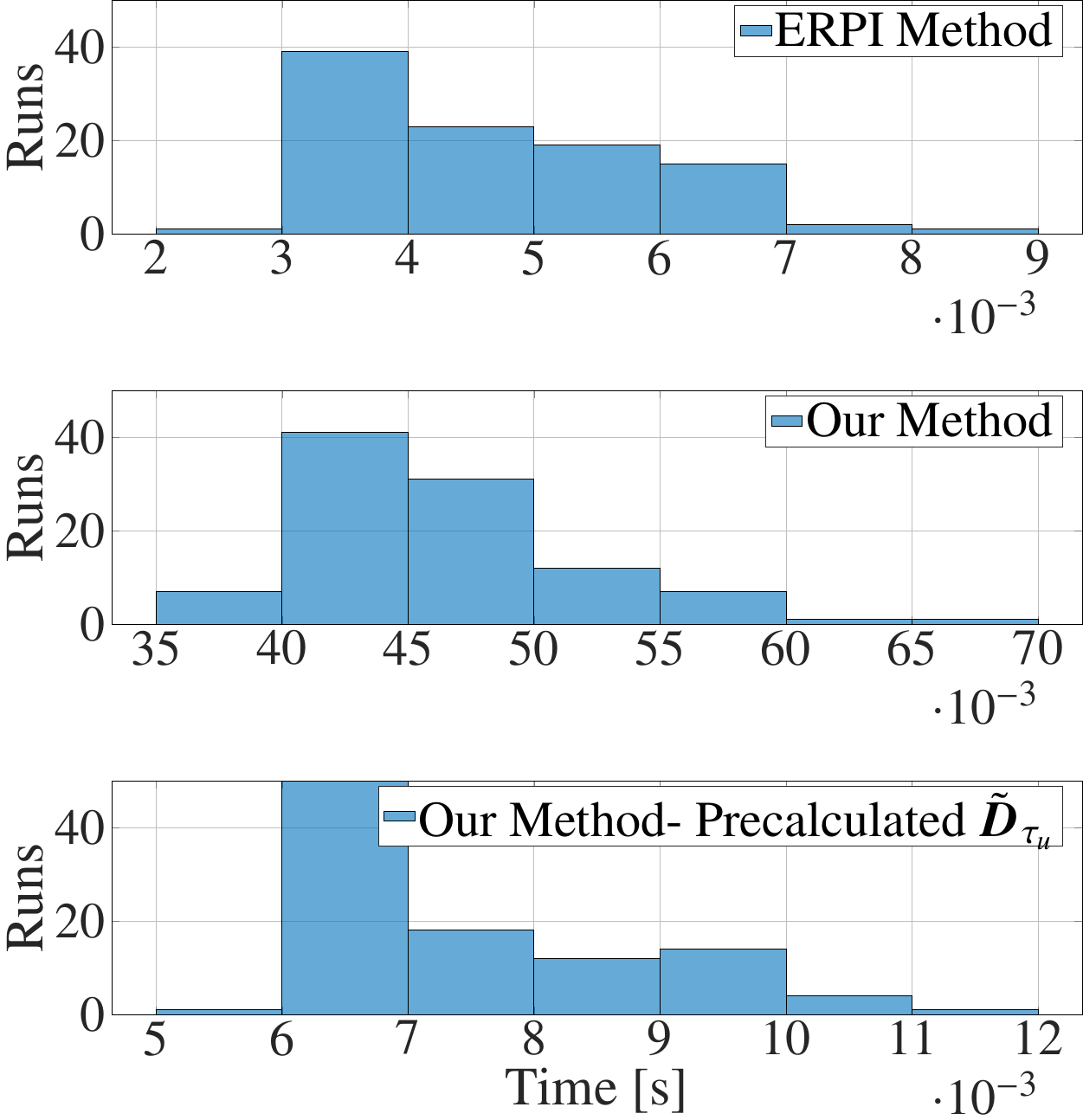}
\caption{Histogram of elapsed run times for the different methods.}
\label{fig:Hist}
\end{subfigure}
\caption{Computational time for the different methods.}
\label{fig:TimeAll}
\end{figure}

 The average elapsed time\footnote{All computations were performed with MATLAB R2022b using an AMD Ryzen 7 PRO 3700U processor with 16 GB of RAM.} over 100 runs of both allocation methods is shown in Fig.~\ref{fig:Time}. Our allocation algorithm is ten times slower than the ERPI algorithm. The calculation of the convex hull and its reformulation into the half-space representation is the main reason for this big discrepancy. In case of constant control effectiveness matrices, this step can be precalculated to decrease the computation time of the allocation, as illustrated in the second bar of Fig.~\ref{fig:Time}. However, given the average elapsed time of 46~\si{ms} for our approach, it can be assumed that real-time capability is still ensured.  Note that the computational time can increase for both methods with an increasing number of actuators. It is also noteworthy that our approach is possibly not programmed efficiently. As illustrated in Fig.~\ref{fig:Hist}, our approach computes 41~\% of the runs within 40--45~\si{ms}, while the ERPI method computes 39~\% of the runs within 3-4~\si{ms}. Precalculation of the convex hull significantly increases the computation time of the allocation problem since 50~\% of the runs are computed within 6--7~\si{ms}. 
\subsection{Real Actuator Dynamics Including Position and Rate Constraints}
\label{sec:SimRes2}
In this section, an analysis is conducted on the performance of the methodology outlined in Sec.~\ref {sec:Method}. The results presented below encompass the real dynamics of the actuator, in addition to the position and rate constraints.

\begin{figure}[H]
\centering
\includegraphics[width=0.7\textwidth]{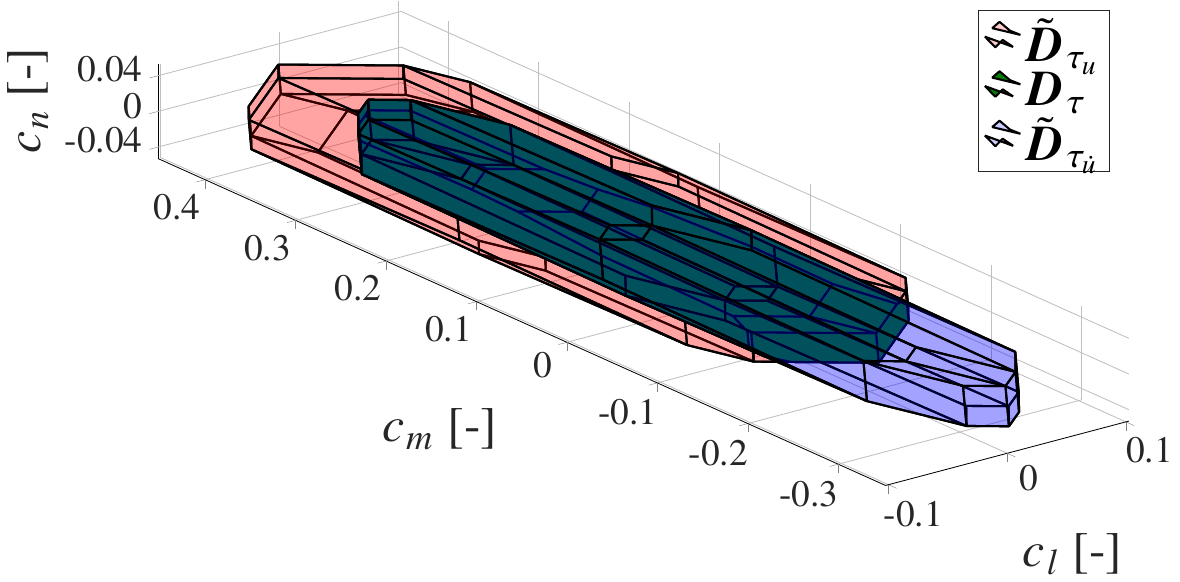}
\caption{AMS $\mat{D}_{\tau}$ from the intersection of two convex hulls.}
\label{fig:AMSNew}
\end{figure}

\begin{figure}[H]
\centering
\includegraphics[width=0.6\textwidth]{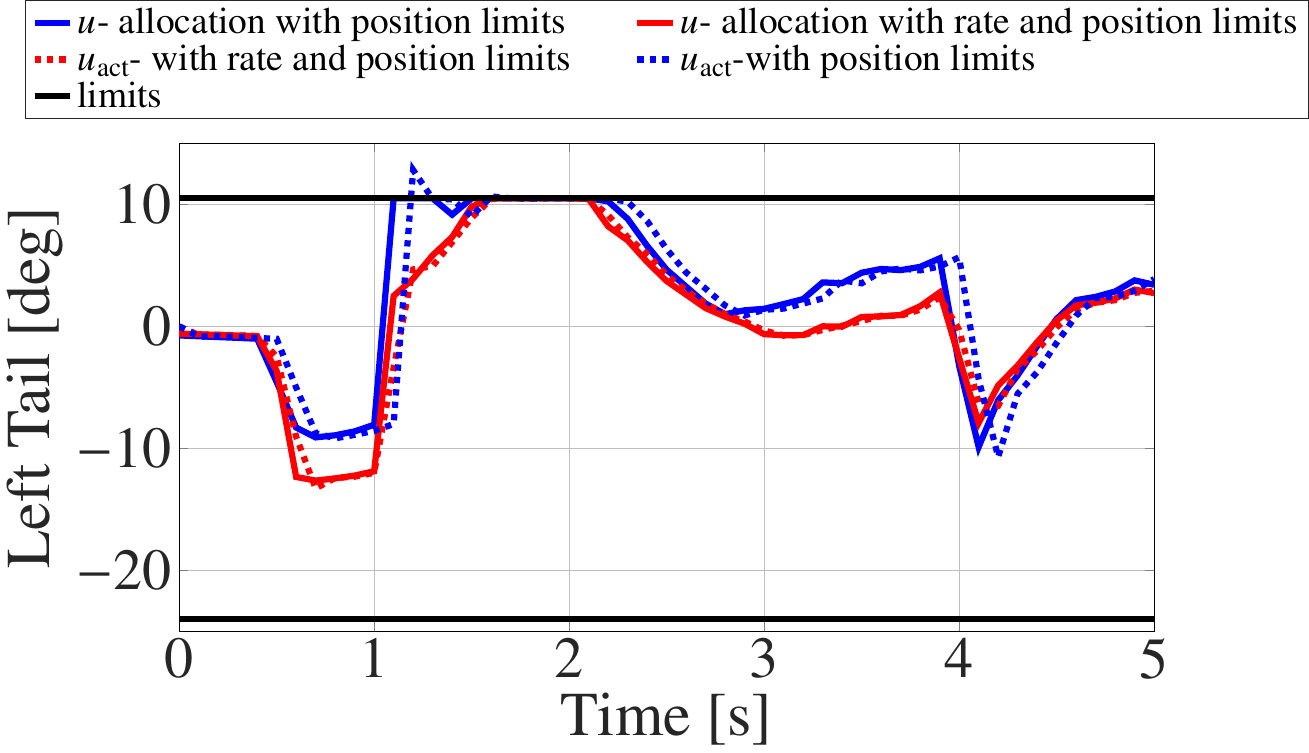}
\caption{Allocation results with and without rate constraints and their realization by the actuator.}
\label{fig:AllocCompare}
\end{figure}

\begin{figure}[H]
    \centering
    \begin{subfigure}[b]{0.45\textwidth}
        \includegraphics[width=\textwidth]{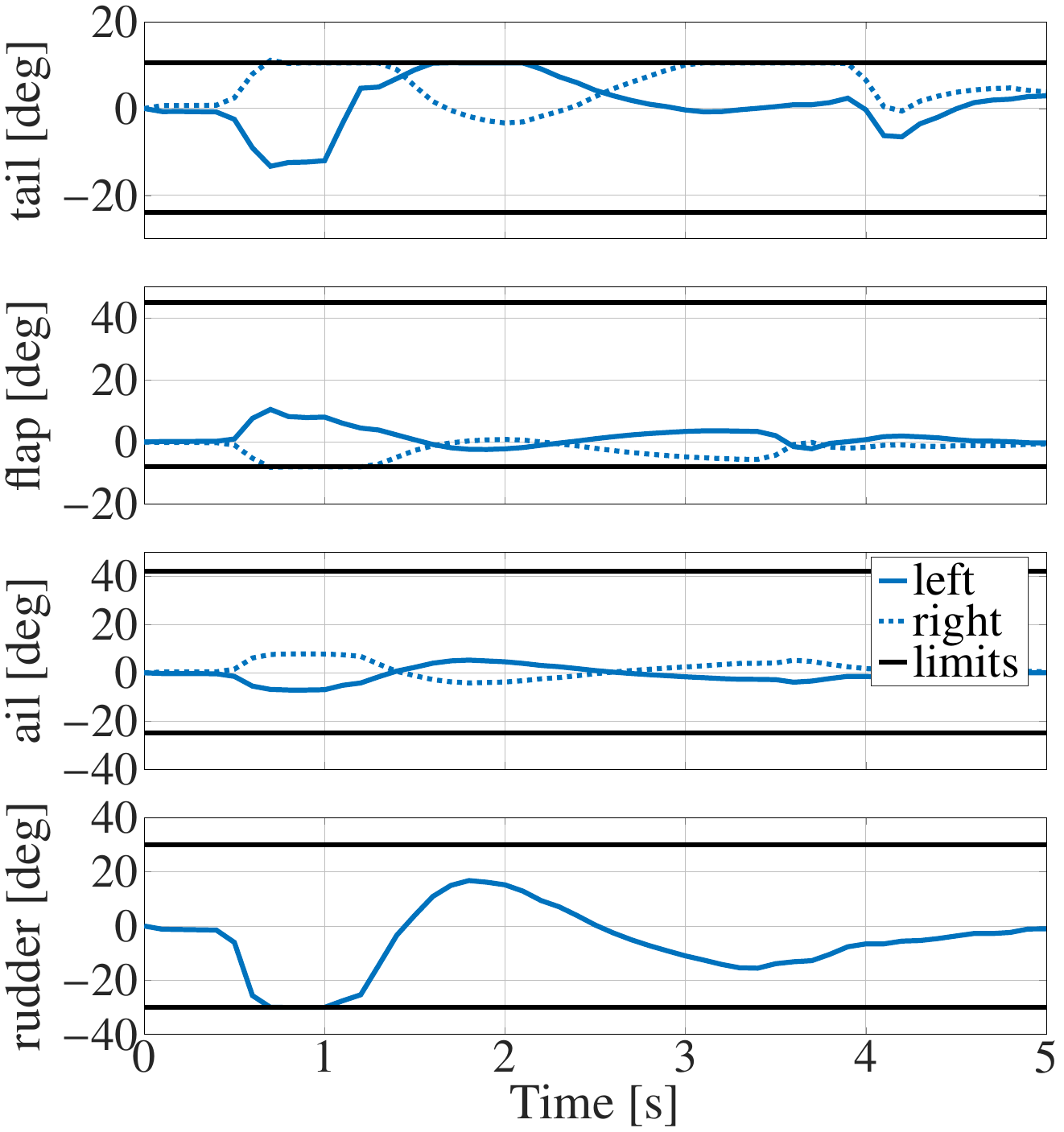}
        \caption{Actuator control input $\vec{u}_{\mathrm{act}}$ with rate and position constraints.}
        \label{fig:RealU}
    \end{subfigure}
    \hfill
    \begin{subfigure}[b]{0.45\textwidth}
        \includegraphics[width=\textwidth]{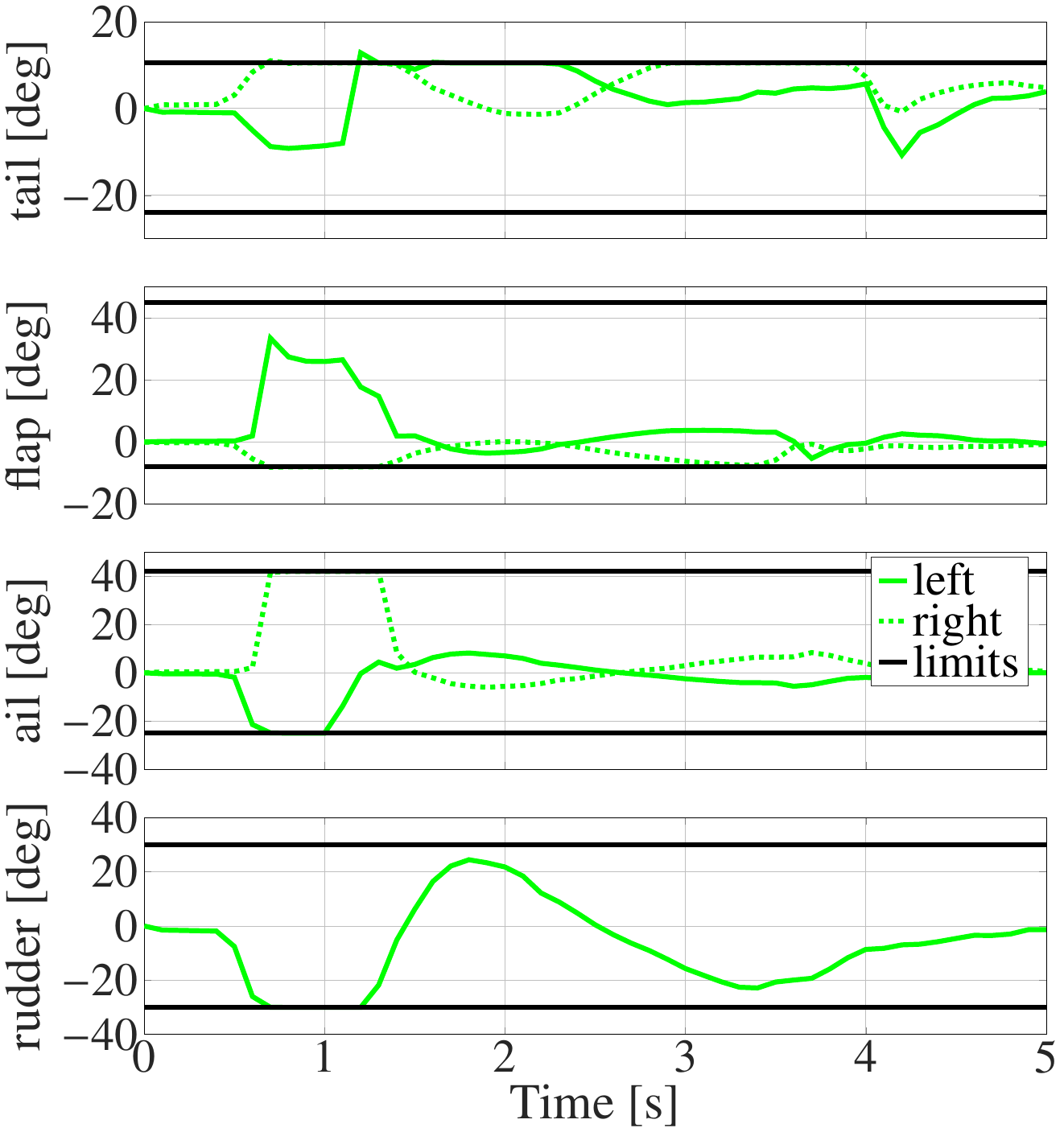}
        \caption{Actuator control input $\vec{u}_{\mathrm{act}}$ with position constraints only.}
        \label{fig:RealUwoRate}
    \end{subfigure}
    \caption{Comparison of the control inputs realized by the actuators.}
    \label{fig:RealUAfterActuator}
\end{figure}

To consider rate constraints, it is necessary to incorporate a first-order homogeneous differential equation into the control allocation algorithm of Eq.~\eqref{eq:ConstrainedControlAlloc}. This differential equation should have a bandwidth that is significantly lower than the bandwidth of the real actuator\footnote{To ensure that the actuator can realize the desired control inputs}. Within this paper, the following state-space form is used to express the differential equation
\begin{equation}
\dot{\vec{u}}=\underbrace{\mathrm{diag}[-2,-2,-2,-2,-2,-2,-2,-2]}_{\mat{A}}\vec{u}. 
\label{eq:Filter}
\end{equation}
Taking into account the matrix $\mat{A}$ from Eq.~\eqref{eq:Filter}, the AMS in Fig.~\ref{fig:AMSNew} is obtained. The AMS in Fig.~\ref{fig:AMSNew} is much smaller than the AMS in Fig.~\ref{fig:AMS}. It is to be noted that the volume of the AMS significantly depends on the choice of $\mat{A}$. Smaller entries in $\mat{A}$ lead to a larger volume of $\tilde{\mat{D}}_{\tau_{\dot{u}}}$ and therefore a larger AMS can be obtained if desired.
To solve the control allocation problem of Eq.~\eqref{eq:ConstrainedAlloc}, \textit{qpoases} is used. The results for the left tail stabilizer in Fig.~\ref{fig:AllocCompare} corroborate the benefit of including rate constraints into the allocation problem.  As demonstrated in Fig.~\ref{fig:AllocCompare}, the solution to the control allocation problem in Eq.~\eqref{eq:ConstrainedControlAlloc} provides the blue line. The determined solution remains within the prescribed position limits. However, the actuator (dashed blue line) fails to adequately execute the control command, resulting in a limit violation and impaired tracking of the desired control command. This is not the case if rate constraints are included in the allocation problem. As shown in Fig.~\ref{fig:AllocCompare}, the actuator precisely follows the red control command and performs more smoothly than the blue line.
The control inputs of the other control surfaces, which are applied by the actuators to the physical system, are depicted in Fig.~\ref{fig:RealUAfterActuator}.
Compared to the sole consideration of position limits in Fig.~\ref{fig:RealUwoRate}, the rate constraints and the dynamical constraint $\dot{\vec{u}}=\mat{A}\vec{u}$ force the actuator to provide smoother control inputs $\vec{u}_{\mathrm{act}}$.
These additional constraints minimize the potential overshoot of the actuators, thereby ensuring precise adherence to the defined position limits. Additionally, the incorporation of dynamic and rate constraints enables precise control inputs with reduced operating time at the aileron and left flap position limitations.

\begin{figure}[H]
    \centering
    \begin{subfigure}[b]{0.45\textwidth}
        \includegraphics[width=\textwidth]{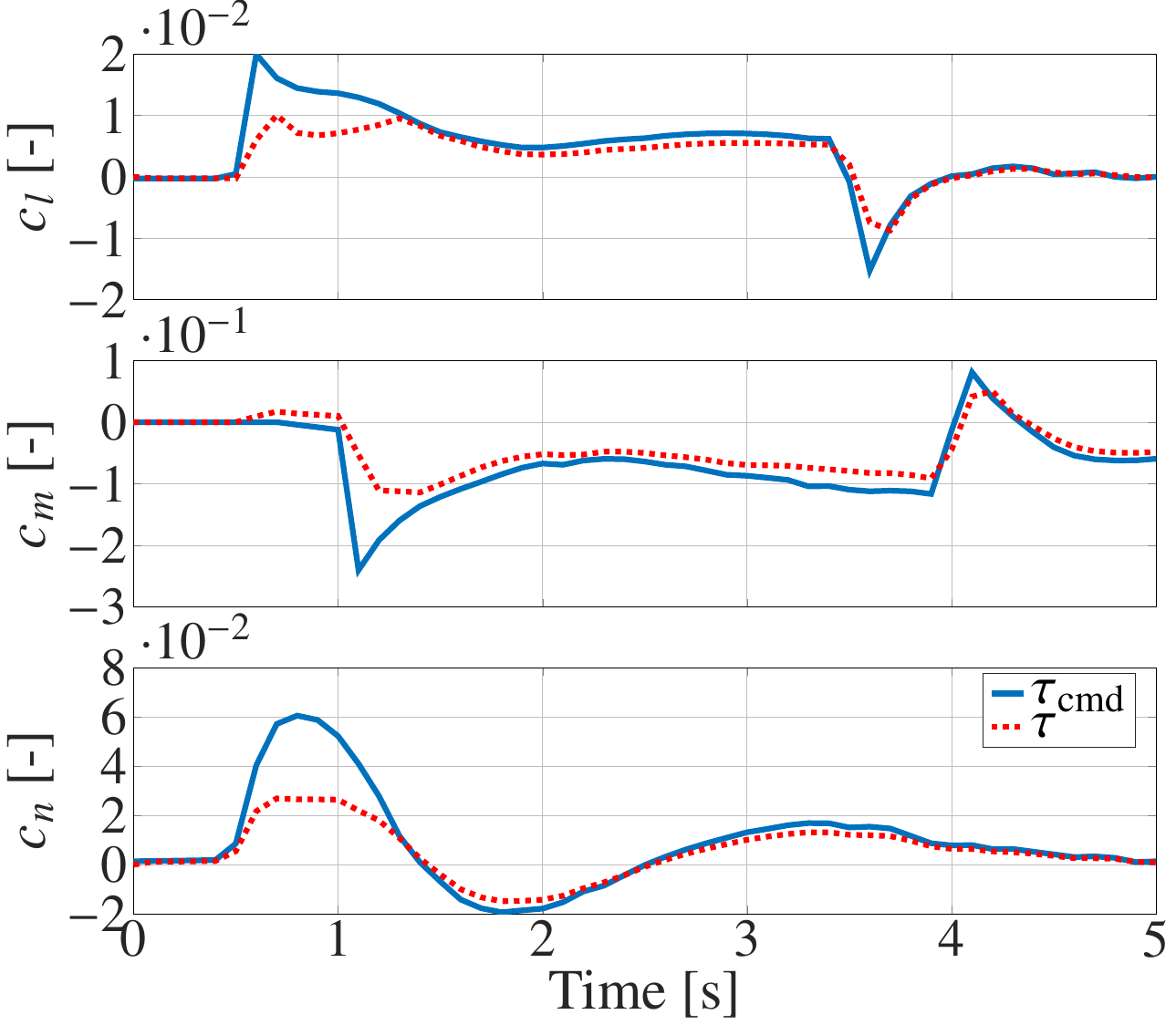}
        \caption{Realized $\vec{\tau}=\mat{B}\vec{u}_{\mathrm{act}}$ with actuator rate and position constraints.}
        \label{fig:RealTauU}
    \end{subfigure}
    \hfill
    \begin{subfigure}[b]{0.45\textwidth}
        \includegraphics[width=\textwidth]{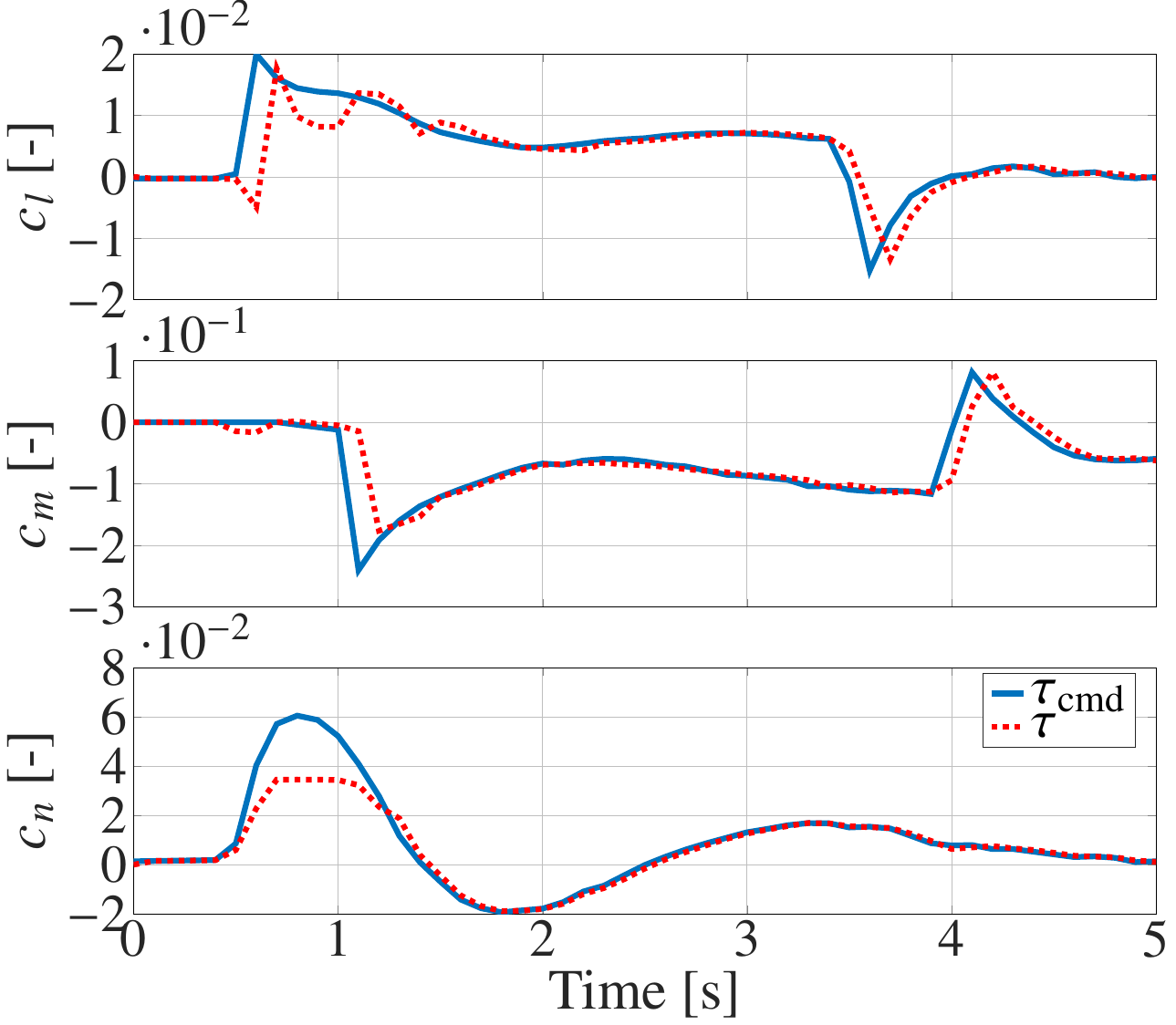}
        \caption{Realized $\vec{\tau}=\mat{B}\vec{u}_{\mathrm{act}}$ with actuator position constraints only.}
        \label{fig:RealTauUwoRate}
    \end{subfigure}
    \caption{Comparison of the implemented control moments.}
    \label{fig:RealTauAfterActuator}
\end{figure}

As depicted in Fig.~\ref{fig:RealTauAfterActuator}, the inclusion of the dynamic and rate constraints into the allocation problem does not significantly affect the control moments. These additional constraints facilitate a smoother progression of the control moments, especially roll moments, compared to a scenario where they are ignored. 
However, a small offset can be seen in Fig.~\ref{fig:RealTauU} for the control moments between two and four seconds. It is conceivable that the method is incapable of accurately tracking the commanded moments without violating any constraints. In general, a comparison of the realized control moments depicted in Fig.~\ref{fig:RealTauAfterActuator} reveals a high degree of similarity, despite the notable disparities in the control inputs.

\section{Conclusion}
\label{sec:Conclusion}
The present paper proposes an alternative QP approach for constrained control allocation that directly incorporates continuous-time rate constraints without the necessity of slack variables. The developed methodology addresses a significant gap in existing control allocation algorithms by ensuring feasibility through the computation of the attainable moment set derived from both position and rate constraints.

A comparison of the proposed methodology with the established ERPI algorithm reveals that they produce equivalent control allocation results for an F18 aircraft while preserving the direction of unattainable commanded moments. This equivalence serves to substantiate the approach's precision and its capacity to accurately replicate the designated attainable moment set. While the proposed method exhibits a computational overhead compared to the ERPI algorithm, the average execution time of 46~\si{ms} may remain within acceptable bounds for real-time applications. For systems with constant control effectiveness matrices, this computational burden can be significantly reduced by performing certain precalculations. Incorporating dynamic and rate constraints into the allocation problem has been demonstrated to result in significant enhancements in actuator behavior. The simulation results demonstrate that these additional constraints prevent operation at actuator limits and lead to smoother control surface movements and enhanced tracking performance. This is essential for ensuring aircraft safety, enhancing ride comfort, and prolonging the lifespan of actuators.

Future research directions include the optimization of the computational implementation to reduce execution time. Moreover, experimental flight tests would provide valuable insights into the practical benefits of this approach in operational environments.
\bibliography{sample}

\end{document}